\documentclass{amsart}
\usepackage{palatino}
\usepackage{euler}

\input xy
\xyoption{all}

\usepackage{latexsym}
\usepackage{amsmath}
\usepackage{amssymb}
\usepackage{amscd}
\usepackage{epsfig}
\usepackage{multicol}
\usepackage{amsfonts}
\usepackage{amsmath}
\usepackage{amsthm}
\usepackage{amssymb}
\usepackage{euscript}
\usepackage{enumerate}
\usepackage{calc}

\newtheorem{THM}{Theorem}

\newenvironment{proofof}[1]{\medskip\noindent
               \textbf{#1.}}{\hfill{{$\square$}}\\}

\newtheorem{theorem}{Theorem}[section]

\newtheorem{cor}[theorem]{Corollary}
\newtheorem{lemma}[theorem]{Lemma}
\newtheorem{prop}[theorem]{Proposition}

\newtheorem{DEF}[theorem]{Definition}

\theoremstyle{definition}
\newtheorem{remark}[theorem]{Remark}

\def\to{\rightarrow}
\def\into{\hookrightarrow}

\numberwithin{equation}{section}

\renewcommand{\O}{\mathcal{O}}

\newcommand{\phitw}{|\!|\Phi |\!|^2|_Q}
\newcommand{\phit}{|\!|\Phi |\!|^2}
\newcommand{\<}{\left<}
\renewcommand{\>}{\right>}

\newcommand{\Tau}{\mathbf{T}}

\newcommand{\T}{T_{\R}}

\newcommand{\Z}{\mathbb Z}

\newcommand{\R}{\mathbb R}
\newcommand{\C}{\mathbb C}

\newcommand{\algt}{\ensuremath{\mathfrak{t}}}

\def\iso{\cong}

\def\t{\mathfrak{t}}

\title{Real loci of symplectic reductions}

\author{R.\ F.\ Goldin}
\address{Mathematical Sciences\\
George Mason University\\ MS 3F2, 4400 University Dr.\\ Fairfax, VA 22030}
\email{rgoldin@math.gmu.edu}

\author{T.\ S.\ Holm}\thanks{The first author was partially supported by NSF grant DMS-0305128. This research was partially conducted
during the period when the second author served as a Clay Mathematics
Institute Liftoff Fellow. The second author was also partially supported by an NSF postdoctoral fellowship}
\address{Department of Mathematics\\University of CA, Berkeley\\ 813 Evans Hall\\ Berkeley, CA 94720}
\email{tsh@math.berkeley.edu}


\begin{document}
\maketitle

\begin{abstract}
Let $M$ be a compact, connected symplectic manifold with a Hamiltonian
action of a compact $n$-dimensional  torus $T$. Suppose that $M$  is
equipped with an  anti-symplectic involution $\sigma$  compatible with
the $T$-action.   The  real  locus of   $M$ is   the fixed point   set
$M^\sigma$   of  $\sigma$.  Duistermaat  introduced   real   loci, and
extended several theorems  of symplectic  geometry  to real  loci.  In
this paper, we extend  another classical result of symplectic geometry
to  real loci:   the  Kirwan surjectivity  theorem.   In addition,  we
compute the kernel of  the real Kirwan  map.  These results are direct
consequences of techniques introduced by Tolman and Weitsman. In some examples, these
results allow us to show that a symplectic reduction $M/\! /T$ has the
same ordinary cohomology as its real locus $(M/\! /T)^{\sigma_{red}}$,
with degrees halved.   This extends  Duistermaat's original result  on
real  loci to a  case in  which there is  not  a natural Hamiltonian torus
action.\\

\end{abstract}

\bibliographystyle{plain}

\section{Introduction}\label{se:intro}

Let $M$ be  a     compact   symplectic  manifold, and     suppose   an
$n$-dimensional torus $T$ acts on $M$ in a Hamiltonian fashion.  Let
$\Phi :M\to\algt^*$ be a moment map.  The components of the moment map
are (equivariant) Morse-Bott functions  on $M$, and  hence may be used
to determine the (equivariant) topology of $M$.

Now suppose that $\sigma$ is an anti-symplectic involution on $M$ that
anti-commutes with the $T$ action:
\begin{equation}\label{eq:anticommute}
\sigma (t\cdot x) =t^{-1}\cdot \sigma(x),
\end{equation}
for all  $t\in  T$ and   $x\in   M$. We assume that   $M^{\sigma}$  is
nonempty.  We call the fixed point set $Q=M^{\sigma}$ of $\sigma$ the
{\em real locus} of  $M$.  The real  locus is a Lagrangian submanifold
of $M$.  The motivating example  of a real  locus is a complex variety
$M$ under complex conjugation $\sigma$.   In this case, the real locus
is the set of real  points on the variety.  By  (\ref{eq:anticommute})
the subgroup of $T$ of elements  of order 2 acts  on $Q$. We call this
subgroup $\T$, the  {\em  real torus}  in  $T$. It  is  immediate that
$\T\iso (\Z/2\Z)^n$.

Duistermaat  \cite{duis} proved that  the  real locus has full  moment
image, that is, that  $\Phi(M)=\Phi(Q)$.  Moreover, he showed  that the
components of the moment map  are Morse-Bott functions for the
real  locus, when using $\Z/2\Z$ coefficients, and  so  we can understand
the  topology of  $Q$ via the 
moment map.  Biss, Guillemin  and the second author \cite{BGH}  proved
that these  moment map components can also  be  used to understand the
equivariant topology of the real locus  with respect to the restricted
$\T$ action.  In addition,  Sjamaar  and O'Shea have  generalized  the
results  of Duistermaat for Hamiltonian    actions of nonabelian  Lie
groups  \cite{OS}.  The principle behind   these results is that  real
loci should behave in a fashion similar to the symplectic manifolds of
which they are submanifolds.  We show  that this philosophy applies in
the context of symplectic reductions.

Given a compact Hamiltonian $T$-space  $M$, Kirwan \cite{kirwan}  proved
that when $T$   acts freely  on  $\Phi^{-1}(\mu)$,  the inclusion  map
$\Phi^{-1}(\mu)\into M$ induces a surjection in equivariant cohomology
with rational coefficients:
\begin{equation}\label{eq:surjection}
\kappa: H_T^*(M)\to H_T^*(\Phi^{-1}(\mu))=H^*(M/\! /T(\mu)).
\end{equation}
The  map $\kappa$ is called the  Kirwan map.  In \cite{TW} the authors
note that under reasonable assumptions about  the torsion of the fixed
point sets  and the  group action,  this map   is surjective over  the
integers as well. For  $T=S^1$ one can  state the result as follows.  For
any prime $p$, assume that {\em either} there is no $p$-torsion in the
$\Z$-cohomology of the fixed point set $M^T$, {\em or} every point not
in     $M^T$    has a     free   $\Z/p\Z$  action.    Then    the  map
(\ref{eq:surjection}) is   surjective  with integer  coefficients. For
higher dimensional $T$, the authors assume  the action is quasi-free to
prove Kirwan's surjectivity over $\Z$, although  this is stronger than
necessary. One  natural question is, what is  the  kernel of $\kappa$?
This question  was answered  in  \cite{TW}  and  refined by  the first
author in the case of rational coefficients \cite{goldin}.

Suppose  that $M$ has the additional  structure  of an anti-symplectic
involution,  compatible with torus  action as specified above. Suppose
also  that $Q:=M^\sigma$ is nonempty.   In  this article  we prove   a
surjectivity    result  for real  loci analogous    to Kirwan's result
(Theorem~1).   We  also use  Tolman and   Weitsman's equivariant Morse
theoretic methods  to compute the kernel  of this  real version of the
Kirwan map (Theorem~2). Thus the main  contribution of this article is
to complete the program begun by  Duistermaat in showing that the real
locus of  a   Hamiltonian $T$-space has  a   $(\Z/2\Z)^n$-action which
behaves as if it were itself a Hamiltonian $T$-space.  Not only is the
equivariant cohomology ring of $Q$ described  by restrictions to fixed
points, but there  is a well-defined  notion of ``real  reduction" and
the major  theorems about  the cohomology  rings of  reduced spaces go
through in the real  case. The proofs  we present are  straightforward
extensions of the work of Atiyah-Bott,  Kirwan, and Tolman-Weitsman to
the real case. We then explore the ramifications of real reduction and
these theorems in a series of examples.

The key to Kirwan's proof of surjectivity is the analysis of the
function
$$
|\! |\Phi |\! |^2:M\to\R.
$$
In  the case of   $S^1$,  this is  not a  Morse  function,  but  it is
Morse-Bott except at  $0$,   the function's minimum  value.  For $\dim
T>1$, there are  a finite number   of critical values of  $\Phi$ where
$|\!   |\Phi|\! |^2$ is    not   Morse-Bott,   in  addition  to    the
minimum. Moreover,
$$
(|\! |\Phi |\! |^2)^{-1}(0)=\Phi^{-1}(0).
$$
Kirwan  uses  the critical  sets  of  $|\! |\Phi   |\! |^2$  to  argue
inductively (on the critical sets)  that the equivariant cohomology of
$M$ surjects onto the ordinary cohomology of $M/\! /T(0)$. It is clear
that surjectivity  for   general values  $\mu$  of $\Phi$   follows by
shifting  the moment map  by an appropriate  constant. These arguments
easily go through when  the function is restricted  to the real locus,
as we will see.

A fundamental step in the proof of surjectivity and the computation of
the kernel is  a lemma due  to  Atiyah and Bott  \cite{AB:local}. This
lemma describes  the local  topology,   in and around  a fixed   point
component   of $M^T$.     Under    the hypothesis  of   no   2-torsion
(Definition~\ref{def:2torsion}) we show  that this  lemma still  holds
for $(\Z/2\Z)^n$-cohomology with   $\Z/2\Z$ coefficients. We  then apply the
lemma and an  inductive argument to both  the question of surjectivity
and to the computation of the kernel for real loci.

For  the  rest  of  this article,  let   $M$  be a compact,  connected
symplectic manifold, endowed  with a Hamiltonian  action  of a compact
torus $T$ and moment map $\Phi$.  Let $\sigma$ be an anti-symplectic
involution, anti-commuting with the torus action, and let $\T\subset T$ be
the set of order 2 elements, plus  the identity. Recall that $\T$ acts
on $Q$.  We begin with a definition.
\begin{DEF}\label{def:2torsion}
Let  $p\in M$ be a  critical point of $\Phi$.  Let $H\subset T$ be the
maximal connected  subtorus of $T$  fixing $p$, $\Lambda_T$ the weight
lattice of $T$, and $\Lambda_H$ the weight lattice of $H$.
Let $N$ be the connected component of   $M^H$  containing $p$.   We
denote  by $\alpha_1,\dots, \alpha_k$ the  weights of the $H$ action
on  $\nu_pN$, the fiber over $p$ of the normal bundle to $N$ in  $M$,
where $k=\dim \nu_pN$. We say that  $p$ is a {\em  2-torsion point} if
$\alpha_i\equiv 0\mod 2\Lambda_H$ for some $i=1,\dots, k$.
\end{DEF}
\begin{remark}
Suppose  $M$ has  a symplectic involution    $\sigma$ with real  locus
$Q$. Suppose $p$ is critical for $\Phi$, fixed by  $H$, and $N$ is the
connected component of  $M^H$ containing $p$. If $p$  lies in $Q$, the
$H$ action on the normal bundle to $N$ in $M$ restricts to an $H_{\R}$
action on the normal bundle $\nu_Q N^\sigma$ to $N^\sigma$ in $Q$. Let
$\beta_1,\dots,\beta_k$   be the    weights     of  the    irreducible
representations obtained  by this action on the  fiber over $p$.  Then
$p$ is a 2-torsion point if and only if  $\beta_i$ is trivial for some
$i$.
\end{remark}

\begin{remark}
 The functions $|\!|\Phi|\!|^2$ and  a  family of perturbations   (see
 Section \ref{se:kernel}) have  critical sets fixed by various subtori
 of  $T$.    Tolman and   Weitsman    prove    $H^*_T(M;\Z)\rightarrow
 H^*(M/\!/T;\Z)$     is    a  surjection    provided    that  $T$ acts
 quasi-freely. As  they note,  one needs  to  ensure that the negative
 normal    bundles to  all critical  sets    of  these functions  have
 nontrivial equivariant Euler classes. For the  real locus case, ``no
 2-torsion points" plays the  same role as  ``quasi-free" does in  the
 symplectic case.  At  the critical sets  of this family of functions,
 the negative normal bundles  have  top Stiefel-Whitney classes  which
 are  nontrivial in  the $\left(H_{\T}^*\otimes 1\right)$-component of
 $H_{\T}^*(N)$  for  each  connected  component  $N$  of the  critical
 set. This amounts to the requirement that, over any point in $N$, the
 negative normal bundle splits  into one-dimensional {\em  nontrivial}
 representations  of   $\T$.  The    condition  that  $Q$ contain   no
 $2$-torsion  points ensures this for   the family of perturbations of
 $|\!|\Phi|\!|^2$.
\end{remark}

The first main theorem is the surjectivity analogue of
\eqref{eq:surjection} for real loci.

\begin{THM}\label{th:surjectivity}
Suppose $M$ is  a compact   symplectic  manifold with a    Hamiltonian
$T^n$-action, and that $Q$ is  the real locus of $M$.  Let $\Phi$ be a
moment map on $M$ and $\mu$ a regular value of $\Phi$. Suppose further
that  $T^n$ acts freely on  $\Phi^{-1}(\mu)$ and  that $Q$ contains no
$2$-torsion     points.   Then   the  {\em     real   Kirwan map}   in
$\T$-equivariant cohomology with $\Z/2\Z$ coefficients
$$
\kappa_{\R}:H^*_{\T}(Q)\to H^*_{\T}(\Phi|_Q^{-1}(\mu))=H^*(Q/\! /\T(\mu)),
$$
induced by inclusion, is a surjection.
\end{THM}

We omit the superscript $n$ in $T^n$ when the dimension is clear.

\begin{remark}
 The  hypothesis that  the real  locus have  no $2$-torsion points  is
 reasonably strong. Real loci of toric  varieties and coadjoint orbits
 in type $A_n$ satisfy this hypothesis, for example, but the real loci
 of maximal coadjoint orbits in type $B_n$ do not.
\end{remark}

The next   task    we  complete  is    to   specify  the     kernel of
$\kappa_{\R}$.  The kernel is precisely the  real analog of the kernel
found in \cite{TW}.  For every $\xi\in\algt$, let
$$
Q_\xi = \left\{ p\in Q\ |\ \left< \Phi(p),\xi\right>\leq 0\right\}.
$$
Let $F=M^T$ denote the fixed point set, and let
$$
K_\xi = \bigg\{ \alpha\in H_{\T}^*(Q;\Z/2\Z)\ \bigg|\ \alpha|_{F\cap
Q_\xi} = 0\bigg\}.
$$
Finally, define the ideal
$$
K_{\R} = \Bigg\langle\sum_{\xi\in\algt}K_\xi\Bigg\rangle.
$$

\begin{THM}\label{th:kernel}
Suppose $M$  is a   compact  symplectic manifold with  a   Hamiltonian
$T$-action, an anti-symplectic involution $\sigma$ anti-commuting
with $T$, and real locus $Q$. Let $\Phi$ be a moment map on $M$ and
$\mu$ a regular value of $\Phi$ such that  $T$ acts freely on
$\Phi^{-1}(\mu)$ and  that $Q$ contains no $2$-torsion  points.
Then the  kernel of  the real Kirwan
map is the ideal $K_{\R}$, so there is a short exact sequence
$$
0\to K_{\R}\to H^*_{\T}(Q)\to H^*(Q/\! /\T(\mu))\to 0.
$$
where the cohomology is taken with $\Z/2\Z$ coefficients.
\end{THM}

The   remainder  of this   paper   is   organized  as  follows.     In
Section~\ref{se:top},    we   review  $G$-equivariant  cohomology with
specific   attention   to  the   case  of   $\Z/2\Z$    coefficients  and
$G=(\Z/2\Z)^n=\T$. Here we make explicit the Thom isomorphism theorem and
the    $\Z/2\Z$      version   of     the    Atiyah-Bott     lemma.    In
Section~\ref{se:morseKirwan}, we study   Morse-Kirwan theory  on  real
loci.  In   Section~\ref{se:surjectivity},  we discuss  reduction,  an
induced anti-symplectic involution on the symplectic reduction, and we
prove the surjectivity theorem (Theorem~\ref{th:surjectivity}).     In
Section~\ref{se:kernel},  we   prove     Theorem~\ref{th:kernel},  the
description  of the  kernel  of  the  real  Kirwan  map.  Finally,  in
Section~\ref{se:EGs},  we  work out   several pertinent  examples.  In
particular, we present an example of a  symplectic reduction which has
the same cohomology as  its real locus, with  degrees divided in half.
As  $M/\!/T$ does not  in general  have a torus  action, this  example
generalizes Duistermaat's original work on real loci, in which he uses
the  $T$ action to  make an analogous statement  for  $M$ and its real
locus.

The authors  would like to  thank Victor Guillemin,  Robert Kleinberg,
Dan Dugger   and Reyer  Sjamaar    for  useful  comments   during  the
preparation of this paper.

\section{Notes on equivariant $\Z/2\Z$-cohomology}\label{se:top}

In  this section we    discuss  equivariant cohomology and  the   Thom
isomorphism in  $\Z/2\Z$  coefficients.  We then state  a  topological
lemma (Lemma   2.1) that  plays an  important  role  in the  proofs of
Theorems 1 and 2.  Given a bundle $E\rightarrow B$ with a $\Z/2\Z$ action
exactly fixing $B$, we show the top  Stiefel-Whitney class $w_k(E)$ is
not a zero-divisor   in  $H_{\Z/2\Z}^*(B;\Z/2\Z)$. Similarly, given   a $\T$
action on a space $E$  that has no  2-torsion, we show that $w_k(E)\in
H_{\T}(B;\Z/2\Z)$ is not a 0-divisor. This  result implies that a certain
long exact  sequence  splits into short   exact sequences (Lemma  2.3)
allowing  inductive arguments to be applied.   This  is not surprising
given the  analogous result due  to Tolman and Weitsman  \cite{TW} for
the $T$-equivariant  cohomology of $B$, where $T$   is a compact torus
acting on $E$ and fixing $B$, and the cohomology is taken with integer
coefficients.

\subsection{Equivariant cohomology with $\Z/2\Z$ coefficients}
Let $G$ be  a compact Lie  group  acting on $M$. Then  the equivariant
cohomology of $M$ is defined to be the ordinary cohomology of the {\em
Borel  construction}. Let $EG$ be  an equivariantly contractible space
on which $G$ acts freely. The Borel construction is
$$
M_G:= M\times_G EG
$$
where the right-hand  side is the set of  equivalence classes of pairs
$(m,e)$ such that $m\in M,\ e\in EG$  and $(m,e)\sim (g\cdot m, g\cdot
e)$ for all $g\in G$. By definition, the $G$-equivariant cohomology of
$M$ is $H_G^*(M):= H^*(M_G)$.  All  $M_G$ are homotopic for  different
choices of  $EG$,  and so any  choice  of $EG$  will  induce  the same
equivariant cohomology ring.

For   any compact  oriented   $G$-manifolds  $N$   and  $M$  and   any
$G$-equivariant  map $f:N\rightarrow  M$,  there is  an associated map
$f_G:N_G\rightarrow M_G$ defined by $f_G(n,e)=(f(n),e)$ on equivalence
classes.   This induces  the   pullback  and    the pushforward  in
cohomology:
\begin{center}
\begin{align*}
f^*&: H_G^*(M)\rightarrow H_G^*(N)\\
f_*&: H_G^*(N)\rightarrow H_G^{*-q}(M)
\end{align*}
\end{center}
where  $q=\dim N   -  \dim M$,  and  may  be  negative. In particular,
equivariant cohomology is functorial. Note  that the pushforward does
not preserve degree, whereas the pullback is a ring map.

In the case of the map $\pi:M\rightarrow pt$ which sends all of $M$ to
a   point $pt$, the  induced   map in  equivariant cohomology  $\pi^*:
H_G^*(p)=H^*(EG/G)\rightarrow H^*_G(M)$ turns $H^*_G(M)$ into a module
over $H^*_G:=H^*_G(pt)$. We note that $H^*_G(M)$  is not always a {\em
free} module over   $H^*_G$: in particular,  if $G$  acts freely  on a
finite dimensional manifold  $M$, then $H^*_G(M)=H^*(M/G)$ (as  can be
seen by noticing that the fiber of $M_G\rightarrow M/G$ is $EG$, which
is contractible) has finite-dimensional cohomology. However, $H^*_G$ is
not finite for nontrivial  $G$, and so the   module is not free.  When
$H^*_G(M)=H^*_G\otimes  H^*(M)$ as vector spaces,  we say  that $M$ is
{\em equivariantly formal}. Clearly,  the  module structure is  free in
this  case.  Hamiltonian   $G$-spaces  are  examples  of equivariantly
formal spaces.

An important fact is that for  $G=(\Z/2\Z)^n$, one choice of $BG$ is $(\R
P^\infty)^n$. Thus, with $\Z/2\Z$ coefficients,
$$
H^*_G= \Z/2\Z[x_1,\dots, x_n]
$$
where $\deg x_i=1$ for all $i$.

\subsection{The Thom isomorphism in $\Z/2\Z$-cohomology}
Let $N\subset M$ be a submanifold of real codimension $k$. In ordinary
cohomology (with any ring of coefficients), the  Thom isomorphism is a
map  between the cohomology    of $N$  and  the  compactly   supported
cohomology of  the  normal bundle $\nu   N$ of $N$  in   $M$. The Thom
isomorphism does not preserve degree; it increases degree by $k$:
$$
\Tau: H^{*}(N)\longrightarrow  H_c^{*+k}(\nu N).
$$
This latter ring is identified using excision to the relative
cohomology as follows:
$$
H_c^*(\nu N) \cong H^*(DN,SN) \cong H^*(M,M-N)
$$
where $DN$  and $SN$ are  the disk and sphere bundles  of  $N$ in $M$,
respectively. The class $\tau=\Tau(1)$  is called the {\em Thom class}
of $\nu N$. Then $\Tau(\alpha)=p^*\alpha\cup \tau$, where $p^*$ is the
map induced by the projection $p:\nu N\rightarrow N$. The class $\tau$
has the  identifying property that, $\int_{\R^k}\tau=1$,  where $\tau$
is restricted to the fiber of $\nu N\rightarrow N$.

In equivariant cohomology  we avoid the  issue  of compact  support by
remaining in the relative cohomology  ring.  The following  discussion
is    an equivariant    version  of    well-known results  found    in
\cite{Spanier} and  \cite{Milnor.Stasheff}. We  now restrict ourselves
to the coefficient ring $\Z/2\Z$.

Let $G$    be  a compact Lie   group  acting  on  $M$. Let    $N$ be a
$G$-invariant submanifold of   $M$  with codimension  $k$  with normal
bundle $\nu N$.   Then  $N_G$ is  a  co-dimension $k$  submanifold  of
$M_G$. Let  $\nu N_G$ be  its normal bundle  and $p:\nu N_G\rightarrow
N_G$ the projection. Note that  $\nu N_G =(\nu  N)_G$ if $\nu N$ is
identified with a $G$-invariant tubular  neighborhood of $N$. Choose a
$G$-invariant Riemannian metric  on   $\nu N_G$,  and let $DN_G$   and
$SN_G$ be the unit disk and sphere bundles, respectively. As $N$ is an
oriented submanifold of $M$, $SN_G$ is  an oriented sphere bundle over
$N_G$. Consider the bundle
$$
\xymatrix{
(D, S) \ar@{^{(}->}[r]^{i}& \ar[d]^{p}  (DN_G, SN_G) \\
& N_G
}
$$
where $(D,S)$ is the fiber over any point $x\in N_G$ in the relative
bundle $(DN_G, SN_G)$.  The cohomology $H^*(DN_G, SN_G)$ is isomorphic
to that of the relative normal bundle $H^*(\nu N_G, \nu N_G-N_G)=:
H^*_G(\nu N, \nu N-N)$.  There is a unique {\em equivariant Thom
class}
$$
\tau_G\in H^k(\nu N_G, \nu N_G-N_G)=: H^k_G(\nu N, \nu N-N)
$$
(in $\Z/2\Z$ coefficients) such that $\tau_G$ restricted to the
(relative) fiber over any point is the unique nonzero class in the
{\em ordinary} cohomology $H^k(D,S)$, where $D$ and $S$ are the
intersection of the fiber with the disk and sphere bundles,
respectively.

In particular, the {\em equivariant Thom isomorphism} in $\Z/2\Z$
coefficients is the map $\Tau: H_G^*(N)\to H_G^{*+k}(\nu N, \nu N-N)$
obtained by the composition of $p^*:H_G^*(N)\to H_G^*(\nu N)$ induced
by the projection, and the map 
$$\cup \tau_G:H_G^*(\nu N)\to
H_G^{*+k}(\nu N, \nu N- N)$$
 which multiplies any class in $H_G^*(\nu
N)$ by $\tau_G$ using the cup product. In other words,
$\Tau(\alpha)=p^*(\alpha)\cup\tau_G$.

\begin{lemma}\label{le:restrictThom}
Let  $N$  be a submanifold of   $M$, where $N$  and   $M$ are compact,
oriented  manifolds.  Let  $\nu N$  be   the normal bundle  of  $N$ in
$M$. Let $G$ act on $M$ fixing $N$, such that for every point $p\in
\nu N-N$, there is a $\Z/2\Z\subset G$ acting nontrivially on $p$.  In
cohomology with $\Z/2\Z$ coefficients, the restriction of the equivariant
Thom class  $\tau_G$  to the  submanifold $N$ is  the top  equivariant
Stiefel-Whitney class of the bundle $\nu N$.
\end{lemma}

\begin{remark}\label{re:2.2}
The reason this  theorem  is so powerful is that it has the  following
interpretation. Let   $f:N\hookrightarrow M$  be  the  inclusion.  The
composition
$$
\xymatrix{
H^*_G(N) \ar[r]^{\hspace{-1in}\Tau} & H^{*+k}_G(\nu N, \nu N-N) \cong
H^{*+k}_G(M, M-N) \ar[r]^{\hspace{1in} j^*} & H^{*+k}_G(M)
}
$$
is the pushforward $f_*$.  Here  $j^*$ is  the  map  induced by  the
inclusion $(M,\emptyset)\subset (M, M-N)$. Lemma~\ref{le:restrictThom}
states that,    in $\Z/2\Z$-cohomology, $f^*f_*(1)=w_k^G(\nu   N)$, where
$f^*$  is the pullback.  In  $\Z$-cohomology, $f^*f_*(1)=e_G(\nu N)$,
the equivariant Euler class of the normal bundle.
\end{remark}

\begin{proof}[Proof of Lemma~\ref{le:restrictThom}]
The  top  Stiefel-Whitney   class  is defined    using  the ``squaring
operation". For $a\in H^k_G(\nu  N, \nu N-N;\Z/2\Z)$, let  $Sq^k(a)=a\cup
a$. The $k^\emph{th}$ {\it equivariant  Stiefel-Whitney class} is  the
unique class $w^G_k(\nu N)\in  H_G^k(N;\Z/2\Z)$ such that $\Tau(w_k^G(\nu
N))= Sq^k(\tau_G)$.  Thus
$$
\Tau(w_k^G(\nu N)) = p^*(w_k^G(\nu N))\cup \tau_G = \tau_G\cup\tau_G.
$$
 Since  $\Tau(f^*j^*(\tau_G))=p^*(f^*j^*(\tau_G))\cup\tau_G=\tau_G\cup
 \tau_G$,    and   $\Tau$ is     an   isomorphism,  it  follows   that
 $f^*j^*(\tau_G)=w_k^G(\nu N)$.
\end{proof}

We   show  in the   proof   of Lemma~\ref{le:atiyahbott}  that,  under
appropriate  assumptions  regarding  the $G$  action on    $\nu N$ (an
assumption implied by no 2-torsion), multiplication by the equivariant
Stiefel-Whitney class is injective.

\subsection{The equivariant top Stiefel-Whitney
class}\label{subsec:2torsion}
Suppose   $M$  is    compact    and  symplectic,   equipped   with  an
anti-symplectic involution  $\sigma$.   We assume   that  there is   a
compact torus $T$ acting  on $M$ in a  Hamiltonian fashion, and a real
torus  $\T  =   (\Z/2\Z)^n\subset    T$   acting  on   the    real  locus
$Q=M^\sigma$. One of the critical properties of Hamiltonian actions on
compact manifolds is the richness of the structure  of the fixed point
sets. We use the behavior of the group action on the normal bundles of
these fixed sets to  obtain  cohomological information. The  following
lemma describes  the cohomology locally near  a fixed point set. It is
crucial in making subsequent inductive arguments.

\begin{lemma}[Atiyah-Bott]\label{le:atiyahbott}
Let  $E\rightarrow N$ be  a rank $k$ bundle  over a  compact, connected, oriented
manifold $N$, and  $G$ either a compact torus  or $G=(\Z/2\Z)^n$
acting  on $E$ whose fixed point set is exactly $N$. Suppose also that
the cohomology of $N$  has no torsion  over $\Z$.  Choose  a
$G$-invariant Riemannian metric  and let   $DE$ and   $SE$  be  the
disk  and  sphere  bundles, respectively, of   $E$.  Then the  long
exact  sequence induced by the inclusions
$(SE,\emptyset)\hookrightarrow (DE, \emptyset)\hookrightarrow (DE, SE)$
splits into short exact sequences
$$
0\rightarrow H^*_G(DE, SE)\rightarrow H^*_G(DE)\rightarrow
H^*_G(SE)\rightarrow 0,
$$
where the coefficient ring is taken to be $\Z$ if $G$ is a compact
torus, and $\Z/2\Z$ for $G=(\Z/2\Z)^n$.
\end{lemma}

\begin{remark}
 In the case  that $G=(\Z/2\Z)^n$, the requirement that   $G$ act on  the
 fibers  of  $E$ nontrivially  is  weaker   than  the condition of  no
 $2$-torsion for $n>1$.
\end{remark}
\begin{remark}
 The condition that   $N$   have no torsion   is  unnecessary  if  the
 coefficient ring is $\Z/2\Z$.
\end{remark}

\begin{proof}
The proof  of \ref{re:2.2} can be found  in \cite{AB:YangMills}  in the  case that
$G=T$, a compact torus.

Suppose that $G=(\Z/2\Z)^n$ and the coefficient field is $\Z/2\Z$. Since $N$ is fixed by $G$ we have
\begin{equation}\label{eq:Ncohom}
H^*_G(N;\Z/2\Z)=H^*_G(pt)\otimes H^*(N) = \Z/2\Z[x_1,\dots,x_n]\otimes H^*(N).
\end{equation}
We show that the leading term of $w_k^G(E)$ in (\ref{eq:Ncohom})
is nonzero, where $\deg (x_i)=1$ and $k$ is the rank of $E$. As we have assumed that $N$ has no
torsion, if the leading term in $H^*_G(pt)\otimes 1$ is not zero, then
the equivariant Stiefel-Whitney class is not a zero-divisor.

Consider the inclusion of a point $p\hookrightarrow N$ and the induced
projection $$ H_G^*(N)\rightarrow H_G^*(p)=H_G^*.$$ We need only show
that the projection of $w_k^G(E)$ is nonzero. The fiber $V$ over $p$
is a representation of $G$ which, since $G$ is simply a product of
$\Z/2\Z$'s, splits into real one-dimensional representations $V_i$ of
$G$.  Each of these representations is nontrivial, by our assumption
that $G$ fixes exactly $N$ in the bundle $E$. The bundle $V_i\times_G
EG\rightarrow BG$ is the pullback of the canonical line bundle over
$\R P^\infty$ under the projection $BG = (\R P^\infty)^n\rightarrow \R
P^\infty$ to one component.  This bundle has nontrivial first
Stiefel-Whitney class by the axiomatic definition of these classes.
Thus the projection of $w_k^G(E)$ is $w_k^G(V) = \prod_i w_1^G(V_i)$
which is nonzero.

We identify $H_T^*(DE,SE)$ with $H_T^{*-k}(N)$ using the Thom isomorphism, and $H_G^*(DE)$ with $H_G^*(N)$ by contraction. Then the map $H_G^*(DE,SE)\rightarrow H_G^*(DE)$ is identified with the injective map given by multiplication by the Euler class (in $\Z$ coefficients) or the
Stiefel-Whitney class (in $\Z/2\Z$ coefficients). It follows that the long
exact sequence
$$
\xymatrix{
\dots\ar[r] &  H^{*-1}_G(SE) \ar[r] & H^*_G(DE, SE) \ar[r] &
H^*_G(DE)\ar[r] & H^*_G(SE)\ar[r] &\dots
}
$$
splits into short exact sequences.

Another proof of this lemma can be found in \cite{allday}.
\end{proof}

\section{Morse-Kirwan theory for real loci}\label{se:morseKirwan}

We return to case where $M$ is a compact connected symplectic manifold
with a Hamiltonian $T$ action with moment map
$\Phi:M\to\algt^*$. Suppose also that $M$ is equipped with an
anti-symplectic involution that anti-commutes with $T$.  Let $Q$ be the
real locus of $M$. We fix an inner product on $\algt^*$.

\subsection{Morse-Kirwan functions on real loci}
In this section, we will prove three important lemmas. In the first, we
determine the critical sets of $\phitw$ in $Q$.  Duistermaat showed
\cite{duis} that for any component of the moment map $\Phi^\xi$,
if $C_M$ is the set of critical points for $\Phi^\xi$ on $M$, the set
of critical points of $\Phi^\xi|_Q$ is $(C_M)^\sigma$.  We prove an
analogous result for $\phitw$.

\begin{lemma}\label{le:critSets}
Let $C_M$ be
the set of critical points of the function $|\!|\Phi|\!|^2:M\to\R$,
and let $C_Q$ be the set of critical points of $|\!|\Phi|_Q|\!|^2:Q\to\R$.  
Then
$C_M$ is preserved by $\sigma$, and $C_M^\sigma=C_Q$.
\end{lemma}

\begin{proof}
Recall that $Q$ is the fixed point set of the involution $\sigma$, and 
that $\sigma$ commutes with the moment map $\Phi$. For any point $p\in Q$ we have $T_pQ=\ker(d\sigma_p-I) = \{v+d\sigma_p(v)|\ v\in T_pM\}$.

Suppose $p$ is critical for $|\!|\Phi|\!|^2$. The involution $\sigma$
induces an isomorphism of tangent spaces.  So $\Phi\circ\sigma = \Phi$
implies that $\sigma(p)$ is critical for $|\!|\Phi|\!|^2$. Thus $C_M$
is preserved by $\sigma$.

Suppose that a $p\in Q$ is critical for $\phit$.  Then
clearly $p$ is critical for $|\!|\Phi|_Q|\!|^2$, so $C_M^\sigma\subseteq 
C_Q$.

Now suppose that $p\in Q$ is critical for $|\!|\Phi|_Q|\!|^2$. Then for 
all $v\in T_pM$, 
\begin{align*}
0&=d|\!|\Phi|_Q|\!|^2_p(v+d\sigma_p(v))\\
&= 2\Phi(p)\cdot d(\Phi|_Q)_p(v+d\sigma_p(v))\\
&=2\Phi(p)\cdot [d(\Phi)_p(v)+d\Phi_p(d\sigma_p(v))]\\
&=2\Phi(p)\cdot [d(\Phi)_p(v)+d(\Phi\circ\sigma)_p(v)]\\
&=2\Phi(p)\cdot [d(\Phi)_p(v)+d(\Phi)_p(v)]\\
&=4\Phi(p)\cdot d\Phi_p(v) = 2d|\!|\Phi|\!|^2_p(v),
\end{align*}
so $p$ is critical for $|\!|\Phi|\!|^2$ as well.
\end{proof}

In the next two lemmas, we show that the functions $\Phi^\xi|_Q$ and
$\phitw$ are Morse-Kirwan functions on $Q$, when $Q$ contains no
$2$-torsion points.   As before, $M$ is
compact, connected, and symplectic with a Hamiltonian torus action and
involution $\sigma$ fixing the real locus $Q$. We suppose also that
$Q$ contains no $2$-torsion points.  The following two lemmas are  real
locus versions of  results of Kirwan \cite{kirwan}. The proof of the
first lemma is nearly identical to that given by Tolman and Weitsman in
\cite{TW:cohom}. The only differences that occur are a result of the
possibility of $2$-torsion points. The proof of the second lemma
contains some additional analysis concerning critical sets. We include
the proofs here for completeness.

\begin{lemma}\label{le:morseKirwan1}
Choose an invariant Riemannian metric on $M$. Given any $\xi\in\algt$,
define $f=\Phi^\xi|_Q$. Let $C$ be a critical set of index $\lambda$
for $f$, and assume that $C$ is the only critical set in the preimage
of an $\varepsilon$-neighborhood around $f(C)$, for some
$\varepsilon>0$.  Define
$$
Q^\pm=f^{-1}\left(\left(-\infty,f(C)\pm \varepsilon\right)\right).
$$
Then the long exact sequence of the pair $\left(Q^+,Q^-\right)$
splits into short exact sequences in equivariant cohomology, with $\Z/2\Z$
coefficients:
$$
0\to H^{*}_{\T}\left(Q^+,Q^-\right) \to
H^*_{\T}\left(Q^+\right) \to H^*_{\T}\left(Q^-\right)\to
0.
$$
\end{lemma}

\begin{proof}
The function $f$ is Morse-Bott at every connected component of the
critical set.  When $\xi$ is generic, the critical sets are
precisely the connected components of the fixed points.  When
$\xi$ is not generic, then the critical sets are connected
components of $(M^K)^\sigma$, where $K\subseteq T$ is some closed subtorus
of $T$.  Since the action is
smooth, these sets are submanifolds.
The negative normal bundle to $C$ is oriented, except at
the minimum, where the negative normal bundle is zero. Consider
the long exact sequence in equivariant cohomology with $\Z/2\Z$
coefficients,
$$
\cdots\to H^{*}_{\T}(Q^+,Q^-)\to H^*_{\T}(Q^+) \to H^*_{\T}(Q^-) \to
H^{*+1}_{\T}(Q^+,Q^-)\to\cdots.
$$
Let $c=f(C)$ be a nonminimal critical value corresponding to $C$.  By
assumption, $c$ is the only critical value in the interval
$[c-\varepsilon,c+\varepsilon]$.  Denote the negative disc and sphere
bundles to $C$ in $Q$ by $D_c$ and $S_c$ respectively.
We let $\lambda$ denote the Morse index of $C$ in $Q$.  Following an
identical argument to \cite[Proof of Proposition 2.1]{TW:cohom},  we
obtain a commutative diagram, with $\Z/2\Z$ coefficients,
\begin{equation}\label{eq:commdiagr}
\begin{array}{c}
\xymatrix{
{}\ar[r] &  H^{*}_{\T}(Q^+,Q^-)\ar[r]^{\gamma_c}\ar[d]_{\iso} &
H^*_{\T}(Q^+) \ar[r]^{\beta_c}\ar[d] & H^*_{\T}(Q^-)\ar[r] &  \\
&  H_{\T}^*(D_c,S_c)\ar[r]^{\delta_c}\ar[d]_{\iso} & H_{\T}^*(D_c) & & \\
&  H^{*-\lambda}_{\T}(D_c) \ar[ur]_{\cup w_c} & & &
}\end{array}
\end{equation}
where $w_c=w_k(D_c)$ is the equivariant top Stiefel-Whitney class of
the bundle $D_c\to C$. The left-most vertical arrows in the diagram are
excision and the Thom isomorphism with $\Z/2\Z$
coefficients. Remark~\ref{re:2.2} ensures that the diagonal arrow is
indeed the cup product with $w_c$. Because there are no $2$-torsion points in
$C\subseteq Q$, $\T$ fixes no sub-bundle  of the negative normal
bundle.  Lemma~\ref{le:atiyahbott} ensures that the cup
product map $\cup w_c$ is injective, and so $\delta_c$ and
$\gamma_c$ must also be injective.  Thus, the long exact sequence
splits.

If $a=f(C)$ is the minimum critical value, the spaces $Q^+$ and $C$
are equivariantly homotopic, and $Q^-$ is empty.  Thus,
$H_{\T}^*(Q^-)=0$ and $H_{\T}^*(Q^+)\iso H_{\T}^*(C)$.  Therefore,
the sequence splits in this case as well.  This completes the proof of
the lemma.
\end{proof}

\begin{lemma}\label{le:morseKirwan2}
Choose an invariant Riemannian metric on $M$. Given $a\in\algt^*$, let
$f_a(x)=\<\Phi(x)-a,\Phi(x)-a\>=|\!|\Phi(x)-a|\!|^2$.
Let $X\subset M$ be a critical set 
for $f_a$, and assume that
$X$ is the only critical set in the preimage of an
$\varepsilon$-neighborhood around $f_a(X)$, for some $\varepsilon>0$.  Define
$$
Q_a^\pm=f_a^{-1}((-\infty,f_a(X)\pm \varepsilon))\cap Q.
$$
Then the long exact sequences of the pair $\left(Q^+_{a},Q^-_{a}\right)$
split into short exact sequences in equivariant cohomology, with $\Z/2\Z$
coefficients:
$$
0\to H^{*}_{\T}\left(Q^+_{a},Q^-_{a}\right) \to
H^*_{\T}\left(Q^+_{a}\right) \to H^*_{\T}\left(Q^-_{a}\right)\to  0.
$$
\end{lemma}

\begin{proof}
Suppose $X$ is not the critical set at the $0$-level set.
While the function $f_a$ may not be Morse-Bott at $X$, Kirwan proves 
\cite{kirwan} that the critical sets behave on a cohomological level as if 
they were nondegenerate. There is a smooth stratification $\{S_\alpha, \alpha\in B\}$ of $M$, where $B$ is an index set of vectors in $\t^*$, such that, for some $\beta\in B$, the stratum $S_\beta$ contains precisely those $x\in M$ such that the limit of the path of steepest descent for $f_a$ from $x$ is in $X$. Kirwan proves that the inclusion $X\subset S_\beta$ is an equivalence of $T$-equivariant cohomology. Furthermore, close to $X$, $S_\beta$ coincides with a $T$-invariant submanifold $\Sigma$ of a neighborhood of $X$, and $\Sigma$ has a well-defined orientable normal bundle.

Let $X^\sigma$ be the fixed point set of the involution $\sigma$ of $X$. All critical sets of $f_a$ restricted to $Q$ will be of the form $X^\sigma$, 
by Lemma~\ref{le:critSets}. Then $X^\sigma$ is the only critical set of $f_a|_Q$ in $Q_a^+$ and not in $Q_a^-$. We assume $X^\sigma$ is connected, for 
each component can be treated individually by a direct sum argument in the 
following discussion. $X^\sigma$ lies in $\Sigma^\sigma$, the submanifold of a neighborhood of $X^\sigma$ in $Q$ formed by the fixed point set of $\sigma$ on $\Sigma$ ($\Sigma^\sigma$ is smooth since $\sigma$ is a smooth action on $\Sigma$). Furthermore, $\Sigma^\sigma$ has a well-defined oriented normal bundle in $Q$, $\Sigma^\sigma$ coincides with $S^\sigma_\beta$ close to $X^\sigma$, and $X^\sigma\subset S^\sigma$ is an equivalence in $T$-equivariant cohomology by \cite{BGH}.

Consider the long exact sequence in equivariant cohomology with $\Z/2\Z$
coefficients,
$$
\cdots\to H^{*}_{\T}(Q_a^+,Q_a^-)\to H^*_{\T}(Q_a^+) \to H^*_{\T}(Q_a^-) \to
H^{*+1}_{\T}(Q_a^+,Q_a^-)\to\cdots.
$$
Except at the minimum, there is a commutative diagram of the
form \eqref{eq:commdiagr}. The $2$-torsion hypothesis guarantees
that $\T$ acts nontrivially on the negative normal bundle to $\Sigma^\sigma$, and
that the cup product with the equivariant Stiefel-Whitney class is an
injection.  By Kirwan's work (extended to the real case as above) the negative normal bundle to $\Sigma^\sigma$ plays the role (cohomologically) that the negative normal bundle to the critical set plays in the case that the critical set is itself a smooth manifold. 

Finally, at the minimum $f_a^{-1}(0)$, the spaces $Q_a^+$
and $X^\sigma$ are equivariantly homotopic, and $Q_a^-$ is empty.  Thus,
$H_{\T}^*(Q_a^-)=0$ and $H_{\T}^*(Q_a^+)\iso H_{\T}^*(X^\sigma)$.
The lemma now follows.
\end{proof}

\subsection{Kirwan's injectivity theorem for real loci}
Let $M$ be a compact symplectic manifold with Hamiltonian $T$-action
with moment map $\Phi$.  Let $M^T$ denote the set of fixed points of
$M$.  Kirwan's injectivity theorem \cite{kirwan} states that the
inclusion
$$
i:M^T\into M
$$
induces an injection in equivariant cohomology
$$
i^*:H_T^*(M)\into H_T^*(M^T).
$$
The real locus analogue is proved in \cite{BGH} using algebraic
techniques. Schmid hints that there is a Morse-Kirwan theoretic proof
in \cite{schmid}, but does not provide the details.

\begin{theorem}[\cite{BGH}, \cite{schmid}]\label{th:Qinjection}
Let $M$ be a compact symplectic manifold with Hamiltonian $T$-action
with moment map $\Phi$.  Let $Q$ be the real locus of $M$ and let
$F=Q^{\T}$ be the fixed points of $\T$ in $Q$, with natural inclusion
map $i:F\into Q$.  Assume that $F$ contains no $2$-torsion points.
Then the pullback in equivariant $\Z/2\Z$-cohomology
$$
i^*:H^*_{\T}(Q)\to H^*_{\T}(F)
$$
is injective.
\end{theorem}

\begin{proof}
Choose an element $\xi\in\algt$ such that the critical set of
$\Phi^\xi$ is precisely $F$.  Order the critical values of
$f=\Phi^\xi$ as $c_1<\cdots<c_N$, and let $F_i$ be the critical points
with $f(F_i)=c_i$. Denote
$$
Q_i^\pm = f^{-1}((-\infty,c_i\pm\varepsilon)).
$$
The theorem holds for $Q_1^-$ because this set is empty.  We proceed
by induction.  Suppose that
$$
H_{\T}^*(Q_i^-)\into H_{\T}^*(Q_i^-\cap F)
$$
is an injection.  We need to apply Lemma~\ref{le:morseKirwan1}, at the
critical sets of $\Phi^\xi$, which are precisely the fixed point
components of $F$.  Thus, since $F$ contains no $2$-torsion points, we
conclude that there is a short exact sequence
$$
0\to H^{*}_{\T}(Q_i^+,Q_i^-)\to H^*_{\T}(Q_i^+) \to H^*_{\T}(Q_i^-) \to 0.
$$

Let $i^*_\pm$ be the inclusion $Q_i^\pm\cap F\into Q_i^\pm$.  Then we
have a commutative diagram, with $\Z/2\Z$ coefficients,
$$
\xymatrix{
0\ar[r] & H^{*}_{\T}(Q_i^+,Q_i^-)\ar[r]\ar[d]^{\iso} & H^*_{\T}(Q_i^+)
\ar[r]\ar[d]^{i^*_+} & H^*_{\T}(Q_i^-) \ar[r]\ar[d]^{i^*_-} & 0 \\
0\ar[r] & H^{*}_{\T}(F_i)\ar[r]^= & H^*_{\T}(Q_i^+\cap F) \ar[r] &
H^*_{\T}(Q_i^-\cap F) \ar[r] & 0.
}
$$
The map $i_-^*$ is an injection by the inductive hypothesis.  A
diagram chase shows that $i_+^*$ is an injection.  Finally, notice
that $Q_i^+$ is equivariantly homotopy equivalent to $Q_{i+1}^-$.
This completes the proof.
\end{proof}

\section{Proving surjectivity (Theorem~\ref{th:surjectivity})}\label{se:surjectivity}

\subsection{Reductions of real loci}
Suppose that $\mu$ is a regular value of $\Phi$, and that $T$ acts
freely on $\Phi^{-1}(\mu)$.  Then $M_{red}=M/\!/T(\mu)$ is a symplectic
manifold with a canonical symplectic form $\omega_{red}$.  Notice that
$\Phi$ is $\sigma$-invariant: $\Phi (\sigma(x))=\Phi(x)$.  Thus,
$\sigma$ acts on the level sets $\Phi^{-1}(\mu)$, taking $T$-orbits to
$T$-orbits.  Therefore, there is an induced involution $\sigma_{red}$
on $M_{red}$. There are natural maps
$$
\xymatrix{
\Phi^{-1}(\mu) \ar@{^{(}->}[r]^{i} \ar[d]^{\pi} & M \\
M/\!/T(\mu) &
}
$$
The symplectic form $\omega_{red}$ on $M/\!/T(\mu)$ satisfies the
property that $\pi^*\omega_{red}=i^*\omega$.  It is immediate that
$\sigma$ also satisfies $\pi^*\sigma_{red}=i^*\sigma$.

\begin{prop}
The involution $\sigma_{red}$ on $M/\!/T(\mu)$ is anti-symplectic with
respect to $\omega_{red}$.
\end{prop}

\begin{proof}
Notice first that
\begin{eqnarray*}
i^*(\sigma^*\omega) & = & -i^*\omega,\ \   \mbox{ since $\sigma^*\omega=-\omega$} \\
                    & = & -\pi^*\omega_{red}  \\
                    & = & \pi^*(-\omega_{red}).
\end{eqnarray*}
But now we also notice that
\begin{eqnarray*}
i^*(\sigma^*\omega) & = & (i^*\sigma)^*(i^*(\omega)),\ \ \mbox{ since $\Phi$ is
                                                        $\sigma$-invariant} \\
                    & = & (\pi^*\sigma_{red})^*(\pi^*(\omega_{red})) \\
                    & = & \pi^*(\sigma_{red}^*\omega^{\phantom{*}}_{red}).
\end{eqnarray*}
Thus $\pi^*(-\omega_{red}) =
\pi^*(\sigma_{red}^*\omega^{\phantom{*}}_{red})$.  Finally, we notice
that $\pi^*$ is an injection, completing the proof.
\end{proof}

\begin{cor}
The real locus $(M/\!/T(\mu))^{\sigma_{red}}$ of $M/\!/T(\mu)$ is a
Lagrangian submanifold.
\end{cor}

\begin{proof}
This follows immediately from the proposition, as Duistermaat has
shown in \cite{duis} that the real locus of any symplectic manifold
with any anti-symplectic involution is Lagrangian.
\end{proof}

The torus $T$ acts freely on $\Phi^{-1}(\mu)$,  so the subset $\T$ of
order two elements acts freely on $(\Phi|_Q)^{-1}(\mu)$. Thus the {\em
reduction} of the real locus is well defined:
$$
Q/\!/\T(\mu):= ((\Phi|_Q)^{-1}(\mu))/\T.
$$

\begin{prop}\label{pr:real}
Let $M$ be a compact symplectic manifold, and suppose $T$ acts on $M$
in a Hamiltonian fashion.  Suppose $M$ has real locus $Q$.  Let
$\mu\in\algt^*$ be a regular value of the moment map such that $T$
acts freely on $\Phi^{-1}(\mu)$.  Then
$$
(M/\!/T(\mu))^{\sigma_{red}}=Q/\!/\T(\mu).
$$
\end{prop}

\begin{proof}
The symplectic reduction $M/\!/T(\mu)$ consists of $T$-orbits with
moment image $\mu$.  The points of $(M/\!/T(\mu))^{\sigma_{red}}$ are
those orbits which are fixed by $\sigma$.  On the other hand, the
reduction $Q/\!/\T(\mu)$ consists of $\T$-orbits in $Q$ with moment
image $\mu$.

Given a $T$-orbit $T\cdot y$ fixed by $\sigma$, there is a $\T$-orbit
in $T\cdot y$ that is fixed pointwise by $\sigma$, for suppose that
$\sigma(y)=t\cdot y$ for some $t\in T$.  Then
$$
x=\sigma(\sqrt{t}\cdot y)=\sqrt{t}\cdot y
$$
is fixed by $\sigma$, and moreover, $\T\cdot x$ is the orbit we seek.
In particular, we can view $\T\cdot x$ as
a $\T$-orbit in $Q$ with moment image $\mu$.

Finally, given a $\T$-orbit $\T\cdot x$ in $Q$ with moment image
$\mu$, this orbit extends to $T\cdot x$, a $T$-orbit fixed by $\sigma$
in $M$.  We claim that this is the only $\T$-orbit extending to
$T\cdot x$.  Suppose that $T\cdot x_1=T\cdot x_2$ for two distinct
$\T$-orbits $\T\cdot x_1$ and $\T\cdot x_2$.  Then there is a $t\in T$
such that $x_1=t\cdot x_2$.  But then, $\sigma(x_1)=\sigma(t\cdot
x_2)$.  Thus, $x_1=t^{-1}\cdot x_2$.  But now, because $T$ acts freely on
$\Phi^{-1}(\mu)$, the element $t$ must be in $\T$, so $\T\cdot
x_1=\T\cdot x_2$, completing the proof.
\end{proof}

\subsection{The proof of surjectivity}

In this section, we will prove that there is a surjection from the
$\T$-equivariant cohomology of $Q$ onto the ordinary cohomology of
$Q/\!/\T(\mu)$.  We will first show this for $\mu=0$, and we can then
deduce the theorem for all other values $\mu$ by shifting the moment
map $\Phi$ appropriately.  We will use the function $f=|\!|\Phi|\!|^2$
as a Morse-Kirwan function on $Q$ and apply Lemma~\ref{le:atiyahbott}
inductively to show that $H_{\T}^*(Q)\to H^*(Q/\!/\T(0))$ with $\Z/2\Z$
coefficients is a surjective map.

Suppose a torus $T$ acts on $M$ in a Hamiltonian fashion. Let $0$ be a
regular value of the moment map $\Phi$, and assume that $T$ acts
freely on $\Phi^{-1}(0)$. Suppose
further that $Q$ is the real locus of $M$ and that $Q$ contains no
$2$-torsion points. Then Theorem~\ref{th:surjectivity} states that the
map with $\Z/2\Z$ coefficients
\begin{equation}\label{eq:realKirwan}
H_{\T}^*(Q)\to H^*(Q/\!/\T(0))
\end{equation}
induced by inclusion is a surjection.  We prove this now.

\begin{proofof}{Proof of Theorem~\ref{th:surjectivity}}
Let $f_0=\phitw$, and let $c_0<c_1<\cdots<c_N$ be the finitely many
critical values of $f_0$, where $c_0=0$. Choose $\varepsilon>0$ so
that
$c_i+\varepsilon<c_{i+1}$ and $c_i-\varepsilon>c_{i-1}$ for all $i$. Let
$$
Q_i^\pm = f_0^{-1}((-\infty,c_i\pm\varepsilon)).
$$
The set $Q_0^-$ is empty, while $Q_0^+$ is equivariantly
homotopy equivalent to $f_0^{-1}(0)$.  The free $\T$ action implies
$$
H_{\T}^*(Q_0^+)=H^*(Q/\!/\T(0)),
$$
in $\Z/2\Z$ coefficients. We note that $Q_i^+$ is equivariantly homotopy equivalent to
$Q_{i+1}^-$ and $Q_N^+=Q$.  As $Q$ has no 2-torsion points, we apply
Lemma~\ref{le:morseKirwan2} to $Q_i^+$ and
$Q_i^-$ to obtain a short exact sequence with $\Z/2\Z$
coefficients
$$
0\to H^{*-\lambda_i}_{\T}(f_0^{-1}(c_i)) \to H^*_{\T}(Q_i^+) \to H^*_{\T}(Q_i^-) \to 0.
$$
Thus, we have a surjection $H^*_{\T}(Q_i^+) \to
H^*_{\T}(Q_{i-1}^+)$ over $\Z/2\Z$, because $Q_i^-$ is equivariantly homotopy
equivalent to $Q_{i-1}^+$.  Over all $i$ we obtain a sequence of surjections in
equivariant cohomology with $\Z/2\Z$ coefficients
$$
H_{\T}^*(Q) = H_{\T}^*(Q_N^+) \to H_{\T}^*(Q_{N-1}^+) \to \cdots
$$
$$
\cdots \to H_{\T}^*(Q_1^+) \to H_{\T}^*(Q_0^+) = H^*(Q/\!/\T(0)),
$$
induced by inclusions $Q_i^+\into Q_{i+1}^+.$
The composition is the Kirwan map (\ref{eq:realKirwan}).
\end{proofof}

We have only proved the theorem for the regular value $\mu=0$.
Suppose $\mu$ is some other regular value of $\Phi$ such that $T$
acts freely on $\Phi^{-1}(\mu)$.  Since $\Phi$ is determined only
up to its constant term, we may choose a new moment map
$\Psi=\Phi-\mu$.  This new moment map will have
$\Psi^{-1}(0)=\Phi^{-1}(\mu)$, and so we may apply the above argument to
$\Psi$ to prove that there is a surjection
$$
\kappa_\R : H_{\T}^*(Q) \to H^*(Q/\!/\T(\mu)).
$$

\section{The kernel of the real Kirwan map (Theorem~\ref{th:kernel})}\label{se:kernel}

We now compute the kernel of the real Kirwan map, $\kappa_{\R}$
(Theorem~\ref{th:kernel}). The technique for real loci is
identical to the symplectic analogue. However, we must impose some
additional requirements on our real locus. First, we require
injectivity of the equivariant cohomology of the real locus to the
equivariant cohomology of its fixed point set
(Theorem~\ref{th:Qinjection}).  Second, we require that the
multiplication by the top Stiefel-Whitney class of the negative
normal bundles to certain functions be nonzero.  This follows
from the proof of Lemma~\ref{le:atiyahbott}. Here the
$(\Z/2\Z)^n$-equivariant Stiefel-Whitney class plays the role that
the equivariant Euler class does in \cite{TW}.

Following \cite{TW} but restricting our attention to the real locus, we define
$$
Q_\xi = \left\{ p\in Q\ |\ \left< \Phi(p),\xi\right>\leq 0\right\}.
$$ for any $\xi\in\algt$, and
$$
K_\xi = \bigg\{ \alpha\in H_{\T}^*(Q;\Z/2\Z)\ \bigg|\ \alpha|_{F\cap Q_\xi} = 0\bigg\},
$$
where $F=Q^{\T}$. Finally, define the ideal
$$
K_{\R} = \Bigg<\sum_{\xi\in\algt}K_\xi\Bigg>.
$$
Theorem~\ref{th:kernel} states that $K_{\R}$ is the kernel of
$\kappa_{\R}$.

\begin{proofof}{Proof of Theorem~\ref{th:kernel}}
First, we show that $K_{\R}\subseteq \ker(\kappa_{\R})$. Let
$\alpha\in K_\xi$ be a class in $K_{\R}$ for some $\xi\in\algt$.
The real-valued function
$f=\Phi^{\xi}|_Q$ is a Morse-Kirwan function on $Q$, in the sense of
Lemma~\ref{le:morseKirwan1}. Let
$d_1<\cdots<d_n$ be the critical values of $f$, and let $i$ be such
that $d_i<0<d_{i+1}$.  Applying Lemma~\ref{le:morseKirwan1} inductively
to the function $f$, we see that $\alpha|_{Q_\xi}=0$.  But
$Q_\xi$ is homotopy equivalent to $Q_{i+1}^-$, where
$Q_{i+1}^-=f^{-1}(-\infty,d_{i+1}-\epsilon)$ for sufficiently small
$\epsilon$. Then $\Phi|_Q^{-1}(0)\subset f^{-1}(0)$ implies
$$
\alpha|_{\Phi|_Q^{-1}(0)}=\alpha|_{f^{-1}(0)}=0.
$$
That is, $\alpha\in\ker(\kappa_{\R})$.

Next, following methods introduced in \cite{TW}, we show that
$\ker(\kappa_{\R})\subseteq K_{\R}$. A class $\alpha\in
\ker(\kappa_{\R})$ implies that $\alpha|_{\Phi^{-1}(0)^\sigma}=0$. By
the injectivity on the real locus (Theorem~\ref{th:Qinjection}), it
suffices to
find $\beta\in K_\R$ such that $\beta|_{F_i}=\alpha|_{F_i}$ for all
components $F_i$ of the fixed point set $Q^{T_\R}$.

Let $f_0=\phitw$, and order the critical sets $C_0,\dots,C_N$ of
$f_0$ with $f_0(C_i)=c_i$ so that the critical values are
$0=c_0<c_1<\cdots<c_N$. Note that $d|\!|\Phi|\!|^2 = d|\!|\
|\!|^2\circ d\Phi$ and so we find the critical sets directly. For any
$x\in Q$ and $v\in T_xQ$ we have
\begin{align}\label{eq:dPhi}
d|\!|\Phi|\!|^2_x(v) &= d|\!|\ |\!|^2_{\Phi(x)}\circ d\Phi_x(v)\\
&= 2\Phi(x)\cdot d\Phi_x(v).
\end{align}
In particular, critical points consist of $x\in \Phi^{-1}(0)$, as well
as any points $x\in Q$ such that $d(\phitw)_x(v)=0$ for all $v\in
T_xQ$, or $d(\phitw)_x(v)$ is perpendicular to $\Phi(x)$ for all $v\in
T_xQ$. Since the differential of $\Phi$ is 0 at fixed points, the
fixed point set  $Q^{\T}$ is critical for $f_0$. Note that
$\alpha|_{C_0}=0$.

Now assume that $\alpha|_{C_i}=0$ for all $i<p$ and $\alpha|_{C_p}\neq
0$.  Applying Lemma~\ref{le:morseKirwan2} to $f_0$, we see that
$\alpha|_{C_p}=mw_k$,
where $mw_k\in H^*_{\T}(C_p)$ is some
multiple of the equivariant top Stiefel-Whitney class $w_k$ of the
negative normal bundle $\nu_{f_0}^-C_p$ to $C_p$ for $f_0$.

Suppose $\alpha|_{C_i}=0$, $i=1,\dots, p-1$. We find $\beta_p\in K_\R$
such that $\alpha|_{C_i}=(\beta_p)|_{C_i}$, $i=1,\dots,p$. Then
$(\alpha-\beta_p)|_{C_i}=0$ for $i=1,\dots,p$ and we apply the
argument inductively to find $\alpha = \sum_{k=p}^N \beta_k\in
K_\R$. Let
$$
\mathcal{S} = \left\{F\in Q_{cc}^{\T}\Big| \ \Phi(F)\cdot
\Phi(C_p)<0\right\},
$$
where $Q_{cc}^{\T}$ denotes the connected components of the fixed
point set.
Note that if $\beta_p\in H^*_{\T}(Q)$ and
$\beta_p|_{\mathcal{S}}=0$, then $\beta_p\in K_\R$.

The class $\beta_p$ is constructed by finding a new Morse-Bott-Kirwan
function $h$ on $Q$ such that:
\begin{enumerate}
\item $Q^{\T}$ is critical for $h$.
\item $C_p$ is critical for $h$, and the negative normal bundle
$\nu_h^-C_p$ equals $\nu_{f_0}^-C_p$.
\item $h(C_i)<h(C_p)$ for all $i=1,\dots, p-1$.
\item $h(F)<h(C_p)$ implies $\Phi(F)\cdot \Phi(C_p)<0$ for any component $F$ of $Q^{\T}$.
\end{enumerate}
Provided such a function $h$ exists, we construct $\beta_p$ as
follows. Let $\overline{C}_0,\dots, \overline{C}_M$ be critical sets
of $h$ such that $h(\overline{C}_0)<\dots<h(\overline{C}_M)$. Since
$C_p$ is critical, we let $(\beta_p)|_{\overline{C}_i}=0$ for all
$\overline{C}_i$ with $h(\overline{C}_i)<h(C_p)$. By
Theorem~\ref{le:morseKirwan2} we may choose
\begin{align*}
\beta_p|_{C_p}&= mw_k(\nu_h^-C_p)\\
&=mw_k(\nu^-_{f_0}C_p)=mw_k \ \mbox{by condition (2)}.
\end{align*}
By condition (1) we have $\beta_p|_F=0$ for all connected components
$F$ in $Q^{\T}$ such that $h(F)<h(C_p).$ It then follows by (3) that
$\beta_p|_{C_i}=0$ for all $i=1,\dots, p-1$ by the argument at the
beginning of this proof. Lastly, we note that by condition (4),
$\beta_p\in K_\R$.

It is left to prove such a function exists.
 This is achieved by perturbing $f_0$ into a family of functions
parametrized by $\lambda\in \R^+$, all of which are equivariantly
perfect, and then we apply Lemma 3.3 to an appropriate perturbation to
find $\beta$. Define
$$f_{\lambda a}=|\!|\Phi|_Q+\lambda a|\!|^2,$$ where $a=c_p$. Note
that for all $x\in Q$ and $v\in T_xQ$,
\begin{align*}
d\left(f_{\lambda a}\right)|_x (v)&= d\left(|\!|\
|\!|^2\right)|_{\Phi|_Q(x)+\lambda a}\circ d(\Phi|_Q+\lambda
a)|_x(v)\\
& =2 \left(\Phi|_Q(x)+\lambda a\right)\cdot d(\Phi|_Q+\lambda a)|_x(v)\\
&=2\left(\Phi|_Q(x)+\lambda a\right)\cdot d(\Phi|_Q)|_x(v)
\end{align*}

Thus $f_{\lambda a}$ is singular on $Q$ exactly when $x\in
\Phi^{-1}(-\lambda a)$, when $d\left(\Phi|_Q\right)|_x(v)=0$ or when
$d\left(\Phi|_Q\right)_x(v)$ is perpendicular to $\Phi|_Q(x)+\lambda
a$. In particular, $f_{\lambda a}$ is singular on $Q^{\T}$, satisfying
condition (1). We note also that $C_p$ is critical for $f_{\lambda a}$
(condition (2)), since for $x\in C_p$ and $v\in T_xQ$,
\begin{align*}
d(f_{\lambda a})|_x(v)&=(\Phi|_Q(x)+\lambda a)\cdot d(\Phi|_Q)_x(v)\\
& =(1+\lambda)a\cdot d(\Phi|_Q)_x(v)\ \mbox{since $\Phi(x)=\Phi(C_p)=c_p=a$}\\
& = (1+\lambda)\frac{1}{2}d(\phitw)|_x(v)\\
&=0, \ \mbox{since $C_p$ is critical for $|\!|\Phi|\!|^2$.}
\end{align*}
This argument also shows that condition (3) holds: the normal bundles $\nu_{f_0}^- C_p$ and
$\nu_{f_{\lambda a}}^-C_p$ are equal, since  $d(f_{\lambda
a})|_x(v)<0$ if and only if $d(f_0)_x(v)<0$ (when $x\in C_p$). Thus
$w_k$ is the equivariant Stiefel Whitney class of the normal bundle to
$C_p$ for $f_{\lambda a}$.  Note that each of the functions in this
perturbed family is still Morse-Bott-Kirwan, in the sense of
Lemma~\ref{le:morseKirwan2}, since in that lemma, we make no
assumption on the value $a$.

Since $f_{\lambda a}(F)=\Phi(F)\cdot \Phi(F) +2\lambda
\Phi(F)\cdot\Phi(C_p)+\Phi(C_p)\cdot\Phi(C_p)$, we note that if
$\Phi(F)\cdot \Phi(C_p)<0$, then for large enough $\lambda$,
$f_{\lambda a}(F) < f_{\lambda a}(C_p)$.
As $Q$ is compact, there is some $L$ such that $\lambda\geq L$ implies
$$
f_{\lambda a}(F) < f_{\lambda a}(C_p)
$$
for all connected components $F$ of $Q^{\T}$ with
$\Phi(F)\cdot \Phi(C_p) <0$.  Thus conditions (1)--(4) are satisfied by $h=f_{La}$. Thus
there exists $\beta\in K_{\R}$ and $\beta_p|_{C_i}=\alpha|_{C_i}$ for
all $i\leq p$, as desired.
\end{proofof}

\section{Examples}\label{se:EGs}

\subsection{The product of two spheres}
Let $\omega$ be the standard symplectic form on $\C P^1$, and let
$S^1$ act on $\C P^1$ by rotation. Let $M=\C P^1\times\C P^1$ with
symplectic form $(\omega,-\omega)$.  Consider the Hamiltonian $T=S^1$
action on $M$ which sends
$$
\theta\cdot (z,w)\mapsto (\theta\cdot z,-\theta\cdot w).
$$
Then the involution $\sigma$ which switches the factors of $M$,
$\sigma (z,w)=(w,z)$ is anti-symplectic and anti-commutes with the
$S^1$ action.  Moreover, the real locus of this involution is
$Q=\{(w,w)\}$, the diagonal copy of $\C P^1$.  The $\T$ action on $Q$
is an action by $\Z/2\Z$, where the nontrivial element acts by rotation
by $\pi$.

Let $\Phi$ be a moment map for the $S^1$ action on $M$, and let $\mu$
be a regular value.  Then $M/\! /S^1(\mu)=\C P^1$ and $Q/\! /\T
(\mu)=\R P^1$.  Theorem~\ref{th:surjectivity} implies that there
is a surjection
$$
\kappa_\R : H_{\Z/2\Z}^*(\C P^1;\Z/2\Z)\to H^*(\R P^1;\Z/2\Z).
$$
The real locus $\C P^1$ is equivariantly formal, since $H^*(\C
P^1;\Z/2\Z)$ has classes only in even degrees.  Thus, as a vector space,
$$
H_{\T}^*(\C P^1;\Z/2\Z)\iso H^*_{\T}\otimes H^*(\C P^1;\Z/2\Z).
$$
The only class in degree $1$ is the equivariant class in
$H^1_{\T}(pt;\Z/2\Z)$, and so the only element of $H_{\T}^*(\C P^1;\Z/2\Z)$
which is not in the kernel of $\kappa_\R$ is this equivariant class.

\subsection{A mod 2 GKM example}

Suppose $M$ is a compact Hamiltonian $T$ space.  Suppose further that $M$ has
only finitely many fixed points, and that for $p\in M^T$, the weights
of the isotropy representation of $T$ on $T_pM$ are pairwise
independent over $\Z$ and over $\Z/2\Z$.  Then we say that
$M$ is a {\em $\!\!\mod 2$ GKM space}.  Such examples include toric
varieties, all coadjoint orbits in type $A_n$, and some (nonmaximal)
coadjoint orbits in other types.  In this case, a result of
\cite{BGH} states that there is an isomorphism
\begin{equation}\label{eq:BGHequality}
H^{2*}_T(M;\Z/2\Z)\iso H^*_{\T}(Q;\Z/2\Z)
\end{equation}
dividing degrees in half.  Moreover, if we write down each of these
rings as subrings of the equivariant $\Z/2\Z$-cohomology of the
appropriate fixed point sets, then the subrings are equal with a
change in grading inherited from the module structure. Both the left-
and right-hand sides of equation (\ref{eq:BGHequality}) are modules
over the corresponding cohomology of a point. In both cases, this
coefficient ring is the polynomial ring $\Z/2\Z[x_1,\dots,x_n]$ but the
generators are regarded as degree 2 classes on the left, and as degree
1 classes on the right.
For example, suppose $M=\O_\lambda$ is a coadjoint orbit of
type $A_2$.
In Figure~1, we show an equivariant class $\alpha$, represented by a
polynomial associated to each fixed point, pictured here on the moment
map image.  As a class in
$H_T^*(\O_\lambda;\Z/2\Z)$, $\alpha$ is a degree $2$ class.  As a class in
$H_{\T}^*((\O_\lambda)^\sigma;\Z/2\Z)$, it is a degree $1$ class.

\begin{figure}[h]
\centerline{
\epsfig{figure=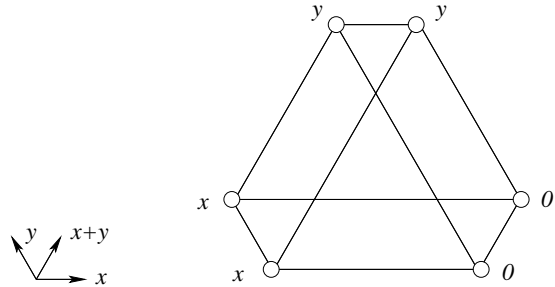,height=1.5in}
}
\smallskip
\centerline{
\parbox{5in}{\caption[coadjoint orbit]{This shows an equivariant
class $\alpha$ on a coadjoint orbit $\O_\lambda$ of type $A_2$.}}
}
\end{figure}

Now let $\Phi:M\to\algt^*$ be the moment map, and suppose
$\mu\in\algt^*$ is a regular value such that $T$ acts freely on
$\Phi^{-1}(\mu)$.  As a result of Kirwan's surjectivity for symplectic
reductions and Theorem~\ref{th:surjectivity}, we have two short exact
sequences in equivariant cohomology, with $\Z/2\Z$ coefficients:
$$
\xymatrix{
0\ar[r] & K \ar[r]^<<<<<{\kappa} & H_T^{2*}(M;\Z/2\Z)\ar[r]\ar@{=}[d] &
    H^{2*}(M/\!/T(\mu);\Z/2\Z)\ar[r] & 0 \\
0\ar[r] & K_{\R} \ar[r]^<<<<<{\kappa_{\R}} & H_{\T}^{*}(Q;\Z/2\Z)\ar[r] &
    H^{*}(Q/\!/\T(\mu);\Z/2\Z)\ar[r] & 0.
}
$$
The Tolman-Weitsman Theorem and Theorem~\ref{th:kernel}
imply that the generators of $K$ and $K_{\R}$ are identified under the equality.
For example, Figure~2 demonstrates the class $\alpha$ which is a
degree $2$ class in $K_\xi(M)$ and is a degree 1 class in $K_\xi(Q)$,
for the $\mu$ and $\xi$ shown.

\begin{figure}[h]
\centerline{
\epsfig{figure=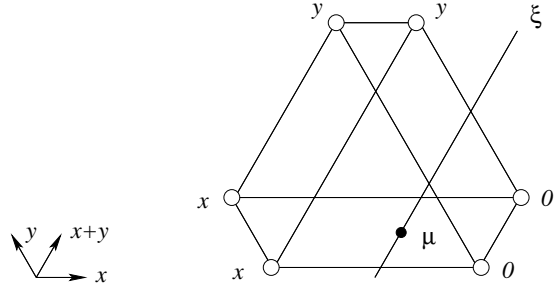,height=1.5in}
}
\smallskip
\centerline{
\parbox{5in}{\caption[coadjoint orbit]{This shows a class
$\alpha$ in $K_\xi$ for $\O_\lambda/\!/T(\mu)$ and for $(\O_\lambda/\!/T(\mu))^\sigma$.}}
}
\end{figure}

Thus, we have a commutative diagram
$$
\xymatrix{
0\ar[r] & K \ar[r]^<<<<<{\kappa}\ar@{=}[d] & H_T^{2*}(M;\Z/2\Z)\ar[r]\ar@{=}[d] &
    H^{2*}(M/\!/T(\mu);\Z/2\Z)\ar[r] & 0 \\
0\ar[r] & K_{\R} \ar[r]^<<<<<{\kappa_{\R}} & H_{\T}^{*}(Q;\Z/2\Z)\ar[r] &
    H^{*}(Q/\!/\T(\mu);\Z/2\Z)\ar[r] & 0.
}
$$
By the Five Lemma,
$$
H^{2*}(M/\!/T(\mu);\Z/2\Z)=H^{*}(Q/\!/\T(\mu);\Z/2\Z)
$$
as rings.  This result is analogous to
Duistermaat's original result on the topology of real loci,; however
$M/\!/T$ may not have a Hamiltonian torus action. Thus this example
extends Duistermaat's result to symplectic reductions.

\subsection{Real toric varieties}
In the symplectic setting, a toric variety is a compact symplectic
$2n$-dimensional manifold with an effective Hamiltonian $T^n$-action.
Let the moment map be $\Phi:M\to \algt^*$.  Then $\Phi(M)=\Delta$ is a
simple rational convex polytope.  Delzant \cite{delzant} proved that
there is a bijection between simple rational convex polytopes in
$\R^n$ and toric varieties.  For every simple rational convex
polytope, Delzant constructed a toric variety $M_\Delta$, which is a
symplectic reduction of affine space.  Namely, $M_\Delta=\C^d/\!/T^k$.

Danilov \cite{danilov} computed the ordinary cohomology of a toric
variety in terms of combinatorial data encoded in its moment polytope.
He showed that
$$
H^*(M_\Delta;\C)=\C[x_1,\dots,x_k]/I,
$$
where $I$ is an ideal of classes depending on how faces of $\Delta$
intersect.  In this case, it is natural to think of the map
$$
H^*_{T^k}(\C^d;\C)=\C[x_1,\dots,x_k]\to H^*(M_\Delta;\C)
$$
as the Kirwan map, and that it is a surjection follows by an argument
similar to Kirwan's.  While $\C^d$ is not compact, the compactness
assumption can be replaced by the assumption that the moment map is
proper and bounded in some direction.  The action of $T^k$ on $\C^d$
certainly satisfies these assumptions.  Indeed, an equivariant
extension of Kirwan's surjectivity, due to the first author
\cite{goldin}, can also be applied, giving a map
$$
H_{T^d}^*(\C^d;\C)\to  H_{T^n}^*(\C^d/\!/T^k;\C).
$$
where $n=d-k$.
The kernel computation does not go through in the same way, but we
have an alternative description of $H_{T^n}^*(M_\Delta;\C)$ which
allows us to compute the kernel for real loci of toric varieties.

We can compute the equivariant cohomology of toric varieties using a
combinatorial description due to Goresky, Kottwitz and MacPherson
\cite{GKM}.  Provided $M_\Delta$ has no 2-torsion, this combinatorial
description also holds for real loci \cite{BGH}, and so we know that
as rings,
$$
H^{2*}_{T^n}(M_\Delta;\Z/2\Z) = H^*_{\T^n}(M_\Delta^\sigma;\Z/2\Z).
$$
Moreover, this equality commutes with the two surjections, giving a commutative diagram
$$
\xymatrix{
0\ar[r] & K=I\ar[r] & H^{2*}_{T^d}(\C^d;\Z/2\Z) \ar@{=}[d]\ar[r] &
                      H^{2*}_{T^n}(\C^d/\!/T^k;\Z/2\Z) \ar@{=}[d]\ar[r] & 0 \\
0\ar[r] & K_{\R}\ar[r] & H^{*}_{\T^d}(\R^d;\Z/2\Z) \ar[r] &
                      H^{*}_{\T^n}(\R^d/\!/\T^k;\Z/2\Z) \ar[r] & 0
          }
$$
Thus, by the Five Lemma, $K=K_{\R}$.  This description of the
cohomology of $(M_\Delta)^\sigma$ as a quotient of a polynomial ring
is identical to Danilov's description of the cohomology of $M_\Delta$.
It was proved using other methods by Davis and Januszkievicz
\cite{DJ}.


\begin{thebibliography}{GKM}

\bibitem{allday} C.\ Allday and V.\ Puppe, {\em Cohomological
Methods in Transformation Groups,} Cambridge Studies in Advanced
Mathematics {\bf 32} 1994, Cambrdige University Press, Cambridge.

\bibitem{atiyah} M.\ Atiyah, Convexity and commuting
Hamiltonians. {\em Bull.\ London Math.\ Soc.} {\bf 14} (1982), no. 1, 1--15.

\bibitem{AB:local}
M.\ Atiyah and R.\ Bott,  The moment map and equivariant
cohomology.  {\em Topology} {\bf 23} (1984), 1--28.

\bibitem{AB:YangMills}
M.\ Atiyah and R.\ Bott, Yang-Mills Equations over Riemann
Surfaces. {\em Philos.\ Trans.\ Roy.\ Soc.\ London A} {\bf 308}, 523-615 (1982).

\bibitem{BV:local}
N.\ Berline and M.\ Vergne,  `Classes caract\'eristiques
\'equivariantes.  Formules de localisation en cohomologie
\'equivariante. {\em C.R.\ Acad.\ Sci.\ Paris S\'er.\ I Math.} {\bf 295}
(1982), 539--541.

\bibitem{BGH}
D.\ Biss, V.\ Guillemin, and T.\ Holm, The mod 2 equivariant
cohomology of fixed point sets of anti-symplectic involutions. To appear in {\em Adv. in Math.}
{\texttt{math.SG/0107151}}.


\bibitem{CS}
T.\ Chang and T.\ Skjelbred, The topological Schur lemma and
related results. {\em Ann.\ of Math. (2)} {\bf 100} (1974), 307--321.

\bibitem{danilov} V.\ Danilov, The geometry of toric varieties. {\em Russian Math.\
Surveys} {\bf 33} (1978), no.\ 2, 97--154.

\bibitem{DJ}
M.\ Davis and T.\ Januszkiewicz.  Convex polytopes, Coxeter
orbifolds and torus actions, {\em Duke Math.\ J.} {\bf 62} (1991), no. 2, 417--451.


\bibitem{delzant} T.\ Delzant, Hamiltoniens p\'eriodiques et images convexes de l'application moment. [Periodic
Hamiltonians and convex images of the momentum mapping]
{\em Bull.\ Soc.\ Math.\ France} {\bf 116} (1988), no.\ 3, 315--339.

\bibitem{duis} H.\ Duistermaat, Convexity and tightness for
restrictions of Hamiltonian functions to fixed point sets of an
anti-symplectic involution. {\em Trans. Amer. Math. Soc.} {\bf 275} (1983), no.\ 1, 417--429.

\bibitem{goldin} R.\ F.\ Goldin,  An effective algorithm for the
cohomology ring of symplectic reductions, {\em Geom.\ and Func.\
Anal.} Vol. 12 (2002), 567--583.



\bibitem{GKM}
M.\ Goresky, R.\ Kottwitz, and R.\ MacPherson,  Equivariant
cohomology, Koszul duality, and the localization theorem,
{\em Invent.\ Math.} {\bf 131} (1998), 25--83.



\bibitem{kirwan}
F.\ Kirwan, {\em Cohomology of Quotients in Symplectic and Algebraic
Geometry}, Princeton University Press, Princeton, NJ, 1984.

\bibitem{Milnor.Stasheff} J.\ Milnor and J.\ Stasheff, {\em
Characteristic Classes}, Princeton University Press, Princeton, NJ,
1974.

\bibitem{OS}
L.\ O'Shea and R.\ Sjamaar, Moment maps and Riemannian symmetric
pairs. {\em Math.\ Ann.} {\bf 317} (2000), no.\ 3, 415--457.

\bibitem{schmid} C.\ Schmid, {\em Cohomologie \'equivariante de
certaines vari\'et\'es hamiltoniennes et de leur partie r\'eelle.}
Th\`ese at Universit\' e de Gen\` eve. \\
 Available at {\texttt
http://www.unige.ch/biblio/these/theses.html}

\bibitem{Spanier} E.\ Spanier, {\em Algebraic Topology}, McGraw-Hill, New York, 1966.

\bibitem{TW}
S.\ Tolman and J.\ Weitsman,  The cohomology ring of
symplectic quotients. {\em Comm. in Anal. and Geom.,} {\bf 11} (2003), no. 4, 751--773.

\bibitem{TW:cohom}
S.\ Tolman and J.\ Weitsman, On the cohomology rings of Hamiltonian
$T$-Spaces. {\em Northern CA Symplectic Geometry Seminar}, 251--258, Amer. Math. Soc. Transl. Ser. 2, vol 196,
Amer. Math. Soc., Providence, RI, 1999.

\end{thebibliography}
\end{document}